\documentclass[10pt]{article}
\usepackage{latexsym}

\newcommand{\duk}{\noindent {\bf Proof. }}
\newcommand{\kduk}{\hfill $\Box$\bigskip}

\newcommand{\N}{\mathbf{N}}

\newtheorem{veta}{Theorem}[section]
\newtheorem{dusl}[veta]{Corollary}
\newtheorem{lema}[veta]{Lemma}
\newtheorem{defi}[veta]{Definition}

\newtheorem{tvrz}[veta]{Proposition}

\def\cla#1#2#3#4#5#6{
  {\sc #1, }#2, {\it #3, }{\bf #4 }(#5), #6.}

\def\kni#1#2#3#4#5{
  {\sc #1, }{\it #2, }#3, #4, #5.}

\def\vsbo#1#2#3#4#5#6#7#8{
  {\sc #1, }#2. In: {#4 (editors), } {\it #5, } #6, #7,  #8; pp. #3.}

\begin{document}

\author{Martin Klazar\thanks{Department of Applied Mathematics (KAM) and Institute for Theoretical 
Computer Science (ITI), Charles University, Malostransk\'e n\'am\v est\'\i\ 25, 118 00 Praha, 
Czech Republic. ITI is supported by the project 1M0021620808 of the 
Ministry of Education of the Czech Republic. E-mail: {\tt klazar@kam.mff.cuni.cz}}}
\title{On growth rates of permutations, set partitions, ordered graphs and other objects}
\date{\today}

\maketitle
\begin{abstract}
For classes ${\cal O}$ of structures on finite linear orders 
(permutations, ordered graphs etc.) endowed with containment order
$\preceq$ (containment of permutations, subgraph relation etc.),
we investigate restrictions on the function $f(n)$ counting objects with size $n$
in a lower ideal in $({\cal O},\preceq)$.
We present a framework of edge $P$-colored 
complete graphs $({\cal C}(P),\preceq)$ which includes many of 
these situations, and we prove 
for it two such restrictions (jumps in growth): $f(n)$ is eventually constant or 
$f(n)\ge n$ for 
all $n\ge 1$; $f(n)\le n^c$ for all $n\ge 1$ for a constant $c>0$ or $f(n)\ge F_n$ for all $n\ge 1$, $F_n$ being the Fibonacci numbers. This generalizes a fragment of a more 
detailed theorem of Balogh, Bollob\'as and Morris \cite{balo-boll-morr06a} on hereditary properties of ordered graphs. 
\end{abstract}

\section{Introduction}

We aim to obtaining general results on jumps in growth of combinatorial structures, 
motivated by such results for permutations \cite{kais_klaz} (which were in turn motivated by results of 
Scheinerman and Zito \cite{schei_zito} and Balogh, Bollob\'as and Weinreich 
\cite{balo-boll-wein00, balo-boll-wein01, balo-boll-wein02} on 
growths of graph properties), and so we begin with them.
{\em Pattern avoidance in permutations}, a quickly developing area of 
combinatorics \cite{atki_murp_rusk, bile, bona_book, bous_stei, clae, eliz_noy, guib_perg_pinz, mans_vain, marc_tard, 
reif, sava_wilf, simi_schm, stan, vatt, west}, is primarily concerned with enumeration of 
sets of permutations 
$$
{\rm Forb}(F)=\{\rho\in{\cal S}:\;\rho\not\succeq\pi\;\forall\pi\in F\},
$$ 
where $F$ is a fixed finite or infinite set of forbidden permutations (patterns) 
and $\preceq$ is the usual containment order on the set of finite permutations 
${\cal S}=\bigcup_{n\ge 0}{\cal S}_n$. Recall 
that $\pi=a_1a_2\dots a_m\preceq\rho=b_1b_2\dots b_n $ iff $\rho$ 
has a subsequence $b_{i_1}b_{i_2}\dots b_{i_m}$, $1\le i_1<i_2<\dots<i_m\le n$, such that 
$a_r<a_s\iff b_{i_r}<b_{i_s}$ for all $1\le r<s\le m$.  

Each set ${\rm Forb}(F)$ is an ideal in $({\cal S},\preceq)$ because 
$\pi\preceq\rho\in {\rm Forb}(F)$ implies $\pi\in {\rm Forb}(F)$ and
each ideal $X$ in $({\cal S},\preceq)$ has the form $X={\rm Forb}(F)$ for some
(finite or infinite) set $F$. For ideals of permutations $X$, it is therefore of interest 
to investigate restrictions on growth of the counting function $n\mapsto |X_n|$, where 
$X_n=X\cap{\cal S}_n$ is the set of permutations with length $n$ lying in $X$. In this direction,   
Kaiser and Klazar \cite{kais_klaz} obtained the following results. 
\begin{enumerate}
\item {\em The constant dichotomy.} Either $|X_n|$ is eventually constant or $|X_n|\ge n$ for all $n\ge 1$. 
\item {\em Polynomial growth.} If $|X_n|$ is bounded by a polynomial in $n$, then there exist integers $c_0,c_1,\dots,c_r$ so that for every $n>n_0$ we have 
$$
|X_n|=\sum_{j=0}^rc_j{n\choose j}.
$$ 
\item {\em The Fibonacci dichotomy.} Either $|X_n|\le n^c$ for all $n\ge 1$ for a constant $c>0$ ($|X_n|$ has 
then the form described in 2) or $|X_n|\ge F_n$ for all $n\ge 1$, where $(F_n)_{n\ge 0}=(1,1,2,3,5,8,13,\dots)$ 
are the Fibonacci numbers.
\item {\em The Fibonacci hierarchy.} The main result of Kaiser and Klazar \cite{kais_klaz} 
states that if $|X_n|<2^{n-1}$ for at 
least one $n\ge 1$ and $X$ is infinite, then there is a unique integer $k\ge 1$ and a constant $c>0$ 
such that 
$$
F_{n,k}\le |X_n|\le n^cF_{n,k}
$$ 
holds for all $n\ge 1$. Here $F_{n,k}$ are the generalized Fibonacci numbers 
defined by $F_{n,k}=0$ for $n<0$, $F_{0,k}=1$, and $F_{n,k}=F_{n-1,k}+F_{n-2,k}+\cdots+F_{n-k,k}$ for $n>0$. 
\end{enumerate}
The dichotomy 3 is subsumed in the hierarchy 4 because $F_{n,1}=1$ and $F_{n,k}\ge F_{n,2}=F_n$ for $k\ge 2$ and $n\ge 1$, but we state 
it apart as it identifies the least superpolynomial growth. Note that the restrictions 
1--4 determine possible growths of ideals of permutations below $2^{n-1}$ but say nothing 
about the growths above $2^{n-1}$. In fact, Klazar \cite{klaz04} showed that while there are only 
countably many ideals of permutations $X$ satisfying $|X_n|<2^{n-1}$ for some (hence, by 4, 
every sufficiently large) $n$, there exists an uncountable family of 
ideals of permutations ${\cal F}$ such that $|X_n|\ll (2.34)^n$ for every $X\in{\cal F}$.

A remarkable generalization of the restrictions 1--4 was achieved by Balogh, Bollob\'as 
and Morris \cite{balo-boll-morr06a} who extended them to ordered graphs. Their main result 
\cite[Theorem 1.1]{balo-boll-morr06a} 
is as follows. Let $X$ be a {\em hereditary property of ordered graphs}, that is, a set 
of finite simple graphs with linearly ordered vertex sets, 
which is closed to the order-preserving graph isomorphism and to the order-preserving
induced subgraph relation. Let $X_n$ be the set of graphs in $X$ with the vertex set 
$[n]=\{1,2,\dots,n\}$. Then, again, the counting function $n\mapsto |X_n|$ is subject to the 
restrictions 1--4 described above. Since ideals of permutations can be represented 
by particular hereditary properties of ordered graphs, this vastly generalizes the results on growth 
of permutations \cite{kais_klaz}.

In this article we present a general framework for proving restrictions of the 
type 1--4 on growths of other classes of structures besides permutations and ordered graphs. 
We shall generalize only 1 and 3, i.e., the constant dichotomy 
(Theorem \ref{bounded}) and the Fibonacci dichotomy (Theorem \ref{fibonacci}). We remark 
that our article overlaps in results with the work of 
Balogh, Bollob\'as and Morris \cite{balo-boll-morr06a}; we explain the overlap presently 
along with summarizing the content of our article. I learned about the results in 
\cite{balo-boll-morr06a} shortly before completing and submitting my work.

We prove in Theorems \ref{bounded} and \ref{fibonacci} that the constant dichotomy and the Fibonacci 
dichotomy hold for ideals of complete graphs having edges colored with $l$ colors, 
where the containment is given by the order-and-color-preserving mappings between vertex
sets. For $l=2$ these structures reduce to graphs with ordered induced subgraph relation and
thus our results on the two dichotomies generalize those of Balogh, Bollob\'as and Morris 
\cite{balo-boll-morr06a} for ordered graphs. To be honest, we must say 
that for the constant dichotomy and the Fibonacci dichotomy it is not hard to reduce the general case $l\ge 2$ to
the case $l=2$ (see Proposition \ref{twocolors} and Corollary \ref{dvebarvy}) and 
so our generalization is not very different from the case of graphs. 
(However, this simple reduction ceases to work for the Fibonacci hierarchy 4.)
Our proofs are different and shorter than the corresponding 
parts of the proof of Theorem 1.1 in \cite{balo-boll-morr06a} (which takes cca 24 pages). 

So instead of (ordered) graphs with induced subgraph relation---which can be captured by complete graphs with edges 
colored in black and white---we consider here complete graphs with 
edges colored in finitely many colors. There is more to this generalization than it might seem, as
we discuss in Section 2, and this is the main contribution of the present article. Our setting enables to capture
many other classes of objects and their containments $({\cal O},\preceq)$ (which need not be directly given in 
graph-theoretical terms) and to show uniformly that their growths are subject to both 
dichotomies. For this one only has to verify (which is usually straightforward) that $({\cal O},\preceq)$ fits the framework of {\em binary classes of objects}. We summarize it briefly now and give details in Section 2.
A binary class of objects is a partial order $({\cal O},\preceq)$ which is realized by embeddings between objects. 
The size of each object $K\in{\cal O}$ is the cardinality of its set of atoms $A(K)$, where an atom of $K$ is an 
embedding of an atom of $({\cal O},\preceq)$ in $K$. For an ideal $X$ in 
$({\cal O},\preceq)$, $X_n$ is the subset of objects in $X$ with size $n$ and we are interested in the counting function 
$n\mapsto|X_n|$. Each set of atoms $A(K)$ carries a linear ordering $\le_K$ and these orderings are preserved 
by the embeddings. The objects $K\in{\cal O}$ and the containment order $\preceq$ are uniquely determined by 
the restrictions of $K$ to the two-element 
subsets of $A(K)$ (the binarity condition in Definition~\ref{binary}). 
Hence $({\cal O},\preceq)$ can be viewed as an ideal in the class $({\cal C}(P),\preceq)$ of 
complete graphs which have edges colored by elements of a finite poset $P$ and where $\preceq$ is the  
edgewise $P$-majorization ordering. For both dichotomies $P$ can be taken without loss of generality to be the 
discrete poset with trivial comparisons. We conclude Section 2 with several examples of binary classes. Here 
we mention three of them. Permutations with the containment of permutations form a binary class. So do finite 
sequences over a finite alphabet $A$ with the subsequence relation. Multigraphs (graphs with possibly repeated 
edges) without isolated vertices and with the ordered subgraph relation form also a binary class; note that their
size is measured by the number of edges rather than vertices.

In Section 3 we prove the constant dichotomy and the Fibonacci
dichotomy for binary classes of objects. In Section 4 we pose some open problems on growths of ideals 
of permutations and graphs and give some concluding comments.

In conclusion let us review some notation. We denote $\N=\{1,2,\dots\}$, $\N_0=\{0,1,2,\dots\}$, 
$[n]=\{1,2,\dots,n\}$ 
for $n\in\N_0$, and
$[m,n]=\{m,m+1,m+2,\dots,n\}$ for integers $0\le m\le n$. For $m>n$ we set $[m,n]=[0]=\emptyset$. If $A,B$ are subsets 
of $\N_0$, $A<B$ means that $x<y$ for all $x\in A,y\in B$. In the case of one-element set we write $x<B$ instead of $\{x\}<B$. 
For a set $X$ and $k\in\N$ we write ${X\choose k}$ for the set of all $k$-element subsets of $X$.  

\bigskip\noindent
{\bf Acknowledgments.} My thanks go to Toufik Mansour and Alek Vainshtein for
their kind invitation to the Workshop on 
Permutation Patterns in Haifa, Izrael in May/June 2005, which gave 
me opportunity to present these results, and to G\'abor Tardos whose insightful remarks 
(he pointed out to me Propositions~\ref{discord} and \ref{twocolors}) helped me to simplify the proofs. 

\section{Binary classes of objects and their examples}

We introduce a general framework for ideals of structures and illustrate it by several examples.

\begin{defi}\label{class}
A {\em class of objects} ${\cal O}$ is given by the following data. 
\begin{enumerate}
\item A countably infinite poset $({\cal O},\preceq)$ that has the least element $0_{\cal O}$. The elements of 
${\cal O}$ are called {\em objects}. We denote the set of atoms of ${\cal O}$ (the objects $K$ such that $L\prec K$ 
implies $L=0_{\cal O}$) by ${\cal O}_1$. ${\cal O}_1$ is assumed to be finite.
\item Finite and mutually disjoint sets $\mathrm{Em}(K,L)$ that are associated 
with every pair of objects $K,L$ and satisfy $|\mathrm{Em}(0_{\cal O},K)|=1$ for every $K$ and 
$\mathrm{Em}(K,L)=\emptyset \iff K\not\preceq L$. The elements of $\mathrm{Em}(K,L)$ are called 
{\em embeddings of $K$ in $L$}.
\item A binary operation $\circ$ on embeddings such that $f\circ g$ is defined whenever $f\in\mathrm{Em}(K,L)$ and 
$g\in\mathrm{Em}(L,M)$ for $K,L,M\in{\cal O}$ and the result is an embedding of $K$ in $M$. This 
operation is associative and has unique left and right neutral elements $\mathrm{id}_K\in\mathrm{Em}(K,K)$. 
It is called a {\em composition of embeddings}.
\item For every object $K\in{\cal O}$ we define
$$
A(K)=\bigcup_{L\in {\cal O}_1}\mathrm{Em}(L,K)
$$
and call the elements of $A(K)$ {\em atoms of} $K$. Each set $A(K)$ is linearly 
ordered by $\le_K$. These linear orders are preserved by the composition: If $f_1,f_2\in A(K)$
and $g\in\mathrm{Em}(K,M)$ for $K,M\in{\cal O}$, then $f_1\le _K f_2 \iff f_1\circ g\le _M f_2\circ g$. 
\end{enumerate}
\end{defi}

\noindent
Note that the set ${\cal O}_1$ is an antichain in $({\cal O},\preceq)$ and that the sets of atoms $A(K)$ 
are finite. To simplify notation, we will use just one symbol 
$\preceq$ to denote containments in different classes of objects.
It follows from the definition that in a class of objects ${\cal O}$ we have $A(0_{\cal O})=\emptyset$ and  
$A(K)=\{\mathrm{id}_K\}$ for every atom $K\in{\cal O}_1$. Every embedding $f\in\mathrm{Em}(K,L)$ induces an 
increasing injection $I_f$ from $(A(K),\le_K)$ to $(A(L),\le_L)$: $I_f(g)=g\circ f$.  For an object $K$ we define its
{\em size} $|K|$ to be the number $|A(K)|$ of its atoms. An {\em ideal} in ${\cal O}$ is any subset $X\subset{\cal O}$
that is a lower ideal in $({\cal O},\preceq)$, i.e., $K\preceq L\in X$ implies $K\in X$. For $n\in\N_0$ we denote
$$
X_n=\{K\in X:\;|K|=|A(K)|=n\}.
$$
Thus $X_0=\{0_{\cal O}\}$. We are interested in the growth rate of the function $n\mapsto|X_n|$ for 
ideals $X$ in ${\cal O}$. 

We postulate the property of binarity.

\begin{defi}\label{binary}
We call a class of objects $({\cal O},\preceq)$ given by Definition~\ref{class} {\em binary} if the following three
conditions are satisfied.
\begin{enumerate}
\item The set ${\cal O}_2=\{K\in{\cal O}:\;|K|=2\}$ of objects with size $2$ is finite. 
\item For any object $K$ and any two-element subset $B\subset A(K)$ the set $R(K,B)=\{L\in{\cal O}_2:\; \exists 
f\in\mathrm{Em}(L,K), I_f(A(L))=B\}$ is nonempty and $(R(K,B),\preceq)$ has the maximum element $M$. We say
that $M$ is the {\em restriction of $K$ to} $B$ and write $M=K|B$.
\item For any object $K$, subset $B\subset A(K)$, and function $h:\;{B\choose 2}\to{\cal O}_2$ such that 
$h(C)\preceq K|C$ for every $C\in {B\choose 2}$, there is a {\em unique} object $L$ with size $|L|=|B|$ such that 
$L|C=h(F(C))$ for every $C\in{A(L)\choose 2}$ where $F$ is the unique increasing bijection from $(A(L),\le_L)$ 
to $(B,\le_K)$. Moreover, for this unique $L$ there is an embedding $f\in\mathrm{Em}(L,K)$ such that $I_f=F$ 
(in particular, $L\preceq K$).
\end{enumerate}
\end{defi}

\noindent
Condition 3 implies that every $K\in{\cal O}$ is uniquely determined by the restrictions to two-element sets of its
atoms (set $B=A(K)$ and $h(C)=K|C$). In particular, in a binary class of objects every set ${\cal O}_n$ is finite. 
If $B\subset A(K)$ and $h(C)=K|C$ for every $C\in{B\choose 2}$, we call the unique $L$ a {\em restriction of $K$ to} $B$ 
and denote it 
$L=K|B$. The full strength of condition 3 for $B\subset A(K)$ and $h(C)\preceq K|C$ is used in the proofs of Propositions~\ref{containment} 
and \ref{universal}.

\begin{tvrz}\label{containment}
In a binary class of objects $({\cal O},\preceq)$, for any two objects $K$ and $L$ we have $K\preceq L$ if and only 
if there is an increasing injection $F$ from $(A(K),\le_K)$ to $(A(L),\le_L)$ satisfying $K|B\preceq L|F(B)$ for every 
$B\in{A(K)\choose 2}$.
\end{tvrz}
\duk
If $K\preceq L$, there exists an $f\in\mathrm{Em}(K,L)$ and by 2 of Definition~\ref{binary} the mapping  
$F=I_f$ has the stated property. In the 
other way, if $F$ is as stated, we define $h:\;{F(A(K))\choose 2}\to{\cal O}_2$ by $h(C)=K|F^{-1}(C)$ and apply 3 of 
Definition~\ref{binary} to $L$, $F(A(K))$, and $h$. The object ensured by it must be equal to $K$ and thus $K\preceq L$. 
\kduk

The main and in fact the only one family of binary classes of objects is given in the following definition.

\begin{defi}\label{edge_colored}
Let $P=(P,\le_P)$ be a finite poset. The class of {\em edge $P$-colored complete graphs} ${\cal C}(P)$ 
is the set of all pairs $(n,\chi)$, where $n\in\N_0$ and $\chi$ is a coloring $\chi:\;{[n]\choose 2}\to P$. The 
containment $({\cal C}(P),\preceq)$ is defined by $(m,\phi)\preceq (n,\chi)$ iff there exists an increasing 
mapping $f:\;[m]\to[n]$ such that for every $1\le i<j\le m$ we have $\phi(\{i,j\})\le_P\chi(\{f(i),f(j)\})$. 
\end{defi}
To show that $({\cal C}(P),\preceq)$ is a binary class of objects one has to specify
what are the embeddings, the composition 
$\circ$, and the linear orders on the sets of atoms, and one has to check that they satisfy 
the conditions in Definitions~\ref{class} and \ref{binary}. This is easy because we modeled  
Definitions~\ref{class} and \ref{binary} to fit $({\cal C}(P),\preceq)$. The least element $0_{{\cal C}(P)}$ is the pair
$(0,\emptyset)$. There is just one atom $(1,\emptyset)$. The embeddings are the increasing mappings $f$ of 
Definition~\ref{edge_colored} and $\circ$ is the usual composition of mappings. If $K=(n,\chi)\in{\cal C}(P)$, it is 
convenient to identify $A(K)$ with $[n]$. Then $\le_K$ is the restriction of the standard ordering of integers. It is 
clear that the conditions of Definition~\ref{class} (properties of embeddings, properties of $\circ$ and the compatibility 
of the orders $\le_K$ and $\circ$) are satisfied. For $K=(n,\chi)\in{\cal C}(P)$ and $B\subset[n]=A(K)$, 
$B=\{a,b\}$ with $a<b$, the restriction $K|B$ is $([2],\psi)$ where $\psi(\{1,2\})=\chi(\{a,b\})$. 
The conditions of Definition~\ref{binary} are easily verified. 

It follows from these definitions that every binary class of objects $({\cal O},\preceq)$ is isomorphic to an ideal  
in some $({\cal C}(P),\preceq)$, up to the trivial distinction that we may have $|{\cal O}_1|>1$ while always 
$|{\cal C}(P)_1|=1$.

\begin{tvrz}\label{universal}
For every binary class of objects $({\cal O},\preceq)$ there is a finite poset $P=(P,\le_P)$ and a 
mapping $F$ from $({\cal O},\preceq)$ to $({\cal C}(P),\preceq)$ with the following properties.
\begin{enumerate}
\item $F$ is size-preserving. 
\item $K\prec L\iff F(K)\prec F(L)$ for every $K,L\in{\cal O}$.
\item $F$ sends all size $1$ objects to $(1,\emptyset)$ but otherwise is injective.
\item $F({\cal O})$ is an ideal in $({\cal C}(P),\preceq)$.
\end{enumerate}
\end{tvrz}
\duk
We set $(P,\le_P)=({\cal O}_2,\preceq)$; $P$ is finite by 1 of Definition~\ref{binary}. If $K\in{\cal O}$ is an object 
with atoms $A(K)=\{a_1,a_2,\dots,a_n\}_{\le_K}$, we 
define $F$ by $F(K)=(n,\chi)$ where $n=|K|$ and, for every $1\le i<j\le n$, $\chi(\{i,j\})=K|\{a_i,a_j\}$. $F$ is 
clearly size-preserving. Also Property 3 is obvious. Property 2 was proved in Proposition~\ref{containment}. 
We prove Property 4. Suppose that 
$(m,\psi)\preceq (n,\chi)=F(K)$ for some $(m,\psi)\in{\cal C}(P)$ and $K\in{\cal O}$. Let $A(K)=\{a_1,a_2,\dots,a_n\}_{\le_K}$.
We take an increasing injection $g:\;[m]\to[n]$ such that $\psi(\{i,j\})\le_P\chi(\{g(i),g(j)\})=K|\{a_{g(i)},a_{g(j)}\}$.
By 3 of Definition~\ref{binary} (applied to $K$, $B=g([m])$, and the $h$ given by $h(C)=\psi(g^{-1}(C))$), there is 
an object $L$, $A(L)=\{b_1,b_2,\dots,b_m\}_{\le_L}$, such that 
$L|\{b_i,b_j\}=\psi(\{i,j\})$ for every $1\le i<j\le m$. Hence $(m,\psi)=F(L)\in F({\cal O})$ and Property 4 is proved.
\kduk

\noindent
Thus ideals in a binary class of objects are de facto ideals in $({\cal C}(P),\preceq)$ for some finite poset
$P$ and it suffices to consider just the classes of objects $({\cal C}(P),\preceq)$. 

The next two results are useful for simplifying proofs of statements on growths of ideals in $({\cal C}(P),\preceq)$.
By a {\em discrete poset} $D_P$ on the set $P$ we understand $(P,=)$, i.e., the poset on $P$ where the only comparisons are equalities. 

\begin{tvrz}\label{discord}
Let $P=(P,\le_P)$ be a finite poset and $D_P$ be the discrete poset on the same set $P$. Then an ideal in 
$({\cal C}(P),\preceq)$ remains an ideal in $({\cal C}(D_P),\preceq)$. 
\end{tvrz}
\duk
Let $X\subset {\cal C}(P)$ be an ideal in $({\cal C}(P),\preceq)$ and let $(m,\psi)\preceq (n,\chi)$ in 
$({\cal C}(D_P),\preceq)$ for some $(m,\psi)\in {\cal C}(P)$ and $(n,\chi)\in X$. By the definitions, then 
$(m,\psi)\preceq (n,\chi)$ in $({\cal C}(P),\preceq)$. So $(m,\psi)\in X$ and $X$ is an ideal in 
$({\cal C}(D_P),\preceq)$ too.
\kduk

\noindent
Thus any general result on ideals in $({\cal C}(D_P),\preceq)$ applies to ideals in $({\cal C}(P),\preceq)$ and 
in many situations it suffices to consider only the simple discrete poset $D_P$.

If $P=(P,\le_P)$ is a finite poset, $b\in P$ is a color, and $D_2=([2],=)$ is the two-element discrete poset, we define a mapping 
$R_b:\;{\cal C}(P)\to {\cal C}(D_2)$ by $R_b((n,\chi))=(n,\psi)$ where $\psi(\{i,j\})=1\iff\chi(\{i,j\})=b$, i.e., we 
recolor edges colored $b$ by $1$ and to all other edges give color $2$. 

\begin{tvrz}\label{twocolors}
Let $X$ be an ideal in $({\cal C}(P),\preceq)$, where $P=(P,\le_P)$ is a finite poset. Then, for every $b\in P$, 
the recolored complete graphs $Y^{(b)}=R_b(X)$ form an ideal in $({\cal C}(D_2),\preceq)$, and for every $n\ge 1$ and 
every color $c\in P$ we have the estimate
$$
|Y^{(c)}_n|\le |X_n|\le\prod_{b\in P}|Y^{(b)}_n|.
$$
\end{tvrz}
\duk
Let $K^*\preceq R_b(L)$ in $({\cal C}(D_2),\preceq)$, where $L\in{\cal C}(P)$. Returning to the original colors, we see
that there is a $K\in{\cal C}(P)$ such that $R_b(K)=K^*$ and 
$K\preceq L$ (even in $({\cal C}(D_P),\preceq)$). This gives the first assertion. The first inequality is 
trivial because the mapping $R_b$ is size-preserving. The second 
inequality follows from the fact that every $K\in{\cal C}(P)$ is uniquely determined by the tuple of values 
$(R_b(K):\;b\in P)$.
\kduk

\noindent
We say that a family ${\cal F}$ of functions from $\N$ to $\N_0$ is {\em product-bounded} if for any 
$k$ functions $f_1,f_2,\dots,f_k$ from ${\cal F}$ there is a function $f$ in ${\cal F}$ such that 
$$
f_1(n)f_2(n)\dots f_k(n)\le f(n)
$$
holds for all $n\ge 1$. Bounded functions, polynomially bounded functions, and exponentially bounded functions are 
all examples of product-bounded families. On the other hand, the family of functions which are, for example, $O(3^n)$ 
is not product-bounded. 

\begin{dusl}\label{dvebarvy}
Let ${\cal F}$ be a product-bounded family of functions and let $g:\;\N\to\N_0$. Suppose that for every ideal 
$X$ in $({\cal C}(D_2),\preceq)$, where $D_2$ is the two-element discrete poset, we have either 
$|X_n|\le f(n)$ for all $n\ge 1$ for some $f\in{\cal F}$ or 
$|X_n|\ge g(n)$ for all $n\ge 1$. Then this dichotomy holds for ideals in every class 
$({\cal C}(P),\preceq)$ for every finite poset $P$.
\end{dusl}
\duk
If $X$ is an ideal in $({\cal C}(P),\preceq)$ and, for $b\in P$, $Y^{(b)}$ 
denotes $R_b(X)$, then either for some $b\in P$ we have 
$|X_n|\ge |Y^{(b)}_n|\ge g(n)$ for all $n\ge 1$ or for every $b\in P$ we have $|Y^{(b)}_n|\le f_b(n)$ for all $n\ge 1$ with
certain functions $f_b\in{\cal F}$. By the assumption on ${\cal F}$ and the inequality in Proposition~\ref{twocolors},  
in the latter case we have $|X_n|\le\prod_{b\in P}f_b(n)\le f(n)$ for all $n\ge 1$ for a function $f\in{\cal F}$.
\kduk

\noindent
We see that to prove for $({\cal C}(P),\preceq)$ an ${\cal F}$-$g$ dichotomy (jump in growth) with a product bounded 
family ${\cal F}$, it suffices to prove it
only in the case $P=D_2$, that is, in the case of graphs with $\preceq$ being the ordered induced subgraph relation. 
This is the case for the slightly weaker version of the constant dichotomy (with $|X_n|\le c$ instead of $|X_n|=c$
for $n>n_0$) and for the Fibonacci dichotomy. On the other hand, the Fibonacci hierarchy, which is an infinite series 
of dichotomies, is a finer result and Corrolary~\ref{dvebarvy} does not apply to it because the corresponding families 
of functions are not product-bounded. 

We conclude this section with several examples of binary classes of objects. 
Our objects are always structures with groundsets $[n]$ for $n$ running 
through $\N_0$ and the containment $\preceq$ is defined by 
the existence of a structure-preserving increasing mapping. Embeddings are these mappings and the 
composition $\circ$ is the usual composition of mappings. With the exception of Examples 7, 8, and 9, the atoms of an object can 
be identified with the elements of its groundset and its size is the cardinality of the groundset. 
We will not repeat these features of $({\cal O},\preceq)$ in every example and we also omit
verifications of the conditions of Definitions~\ref{class} and \ref{binary} which are easy. With the exception of Example 6, 
each set $R(K,B)$ of 2 of Definition~\ref{binary} has only one element and condition 2 is satisfied automatically. In every example
we mention what is the poset $(P,\le_P)=({\cal O}_2,\preceq)$ (see Proposition~\ref{universal}). It is  
the discrete ordering $D_k=([k],=)$ for some $k$, with exception of Example 6 where it is the linear ordering $L_2=([2],\le)$. 
In Example 6 the sets $R(K,B)$ have one or two elements. In Examples 7, 8, and 9 the atoms are edges rather than vertices and 
the size of an object is the number of its edges. 

\bigskip\noindent
{\bf Example 1. Permutations.} ${\cal O}$ is the 
set of all finite permutations, which are the bijections $\rho:\;[n]\to[n]$ where $n\in\N_0$. For two permutations 
$\pi:\;[m]\to[m]$ and $\rho:\;[n]\to[n]$, we define $\pi\preceq\rho$ iff there is an increasing mapping $f:\;[m]\to[n]$ such 
that $\pi(i)<\pi(j)\iff \rho(f(i))<\rho(f(j))$; this is just a reformulation of the definition given in the beginning of Section 1.
There is only one atom, the $1$-permutation, and ${\cal O}_2$ consists of the two $2$-permutations. $(P,\le_P)$ is the 
discrete ordering $D_2$. By Proposition~\ref{universal}, permutations form an ideal in $({\cal C}(D_2),\preceq)$. 
It is defined by the ordered transitivity of both colors: if $x<y<z$ and $\{x,y\}$ and $\{y,z\}$ are colored $c\in[2]$, 
then $\{x,z\}$ is colored $c$ as well.

\bigskip\noindent
{\bf Example 2. Signed permutations.} We enrich permutations $\rho:\;[n]\to[n]$ by coloring the elements of the definition domain 
$[n]$ white ($+$) and black ($-$), and we require that the embeddings $f$ are in addition color-preserving. There are two atoms and
${\cal O}_2$ consists of eight signed $2$-permutations. $(P,\le_P)$ is the discrete ordering $D_8$.

\bigskip\noindent
{\bf Example 3. Ordered words.} ${\cal O}$ consists of all mappings $q:\;[n]\to[n]$ such that the image of 
$q$ is $[m]$ for some $m\le n$. For two such mappings $p:\;[m]\to[m]$ and $q:\;[n]\to[n]$ we define $p\preceq q$ in the 
same way as for permutations. The elements of $({\cal O},\preceq)$ can be viewed as words  
$u=b_1b_2\dots b_n$ such that 
$\{b_1,b_2,\dots,b_n\}=[m]$ for some $m\le n$, and $u\preceq v$ means that $v$ has a subsequence with the same length as $u$ 
whose entries form the same pattern (with respect to $<,>,=$) as $u$. There is one atom and 
${\cal O}_2$ consists of three elements ($12$, $21$, and $11$). $(P,\le_P)$  is the discrete ordering $D_3$. 

\bigskip\noindent
{\bf Example 4. Set partitions.} ${\cal O}$ consists of all partitions $([n],\sim)$ where $\sim$ is an equivalence 
relation on $[n]$. We set $([m],\sim_1)\preceq([n],\sim_2)$ iff there is a subset $B=\{b_1,b_2,\dots,b_m\}_<$ of 
$[n]$ such that $b_i\sim_2 b_j\iff i\sim_1 j$. There is only one atom and ${\cal O}_2$ 
has two elements. $(P,\le_P)$  is the discrete ordering $D_2$. By Proposition~\ref{universal}, partitions form an ideal 
in $({\cal C}(D_2),\preceq)$. It is defined by the transitivity of the color $c$ corresponding to the partition 
of $[2]$ with $1$ and $2$ 
in one block: If $x,y,z$ are three distinct elements of $[n]$ such that $\{x,y\}$ and $\{y,z\}$ are colored $c$, 
then $\{x,z\}$ is colored $c$ as well. To put it differently, set partitions can be represented by ordered graphs whose components are complete graphs. Pattern avoidance in set partitions was investigated by Klazar \cite{klaz00}, for futher results see Goyt \cite{goyt} and Sagan \cite{saga}.

\bigskip\noindent
{\bf Example 5. Ordered induced subgraph relation.} ${\cal O}$ is the set of all simple graphs with vertex 
set $[n]$. For two graphs $G_1=([n_1],E_1)$ and $G_2=([n_2],E_2)$ we define $G_1\preceq G_2$ iff there is an 
increasing mappings $f:\;[n_1]\to[n_2]$ such that $\{x,y\}\in E_1\iff \{f(x),f(y)\}\in E_2$. Thus $\preceq$ is the ordered 
induced subgraph relation. There is only one atom and ${\cal O}_2$ has two elements. $(P,\le_P)$ is is the discrete ordering 
$D_2$. This class essentially coincides with $({\cal C}(D_2),\preceq)$.

\bigskip\noindent
{\bf Example 6. Ordered subgraph relation.} We take ${\cal O}$ as in the previous example and in the definition 
of $\preceq$ we change $\iff$ to $\Longrightarrow$. Thus $\preceq$ is the ordered subgraph relation. 
There is only one atom and ${\cal O}_2$ 
has two elements. Unlike in other examples, $({\cal O}_2,\preceq)$ is not a discrete ordering but the linear ordering 
$L_2$. Every set $R(K,B)$, where $K$ is a graph and $B$ is a two-element set of its vertices (atoms), has one or two elements 
and $(R(K,B),\preceq)$ is $L_1$ or $L_2$. Thus $(P,\le_P)$ is the linear ordering $L_2$. This class essentially coincides with 
$({\cal C}(L_2),\preceq)$.

\bigskip\noindent
{\bf Example 7. Ordered graphs counted by edges.} Let ${\cal O}$ be the set of simple graphs with the vertex set $[n]$ and 
without isolated vertices, and let $\preceq$ be the ordered subgraph relation (as in the previous example). There is one atom 
corresponding to the single edge graph. 
The size of $G=([n],E)$ is now $|E|$, the number of edges. ${\cal O}_2$ has six elements and $({\cal O}_2,\preceq)$ is 
$D_6$. The linear ordering $\le_G$ on $E$, 
the set of atoms of $G=([n],E)$, is the restriction of the lexicographic ordering $\le_l$ on ${\N\choose 2}$: 
$e_1\le_l e_2\iff\min e_1<\min e_2$ or ($\min e_1=\min e_2\;\&\;\max e_1<\max e_2$). It is clear that $\le_l$ is compatible with 
the embeddings, which are increasing mappings between vertex sets sending edges to edges, and so condition 4 of 
Definition~\ref{class} is satisfied. Let us check the conditions of Definition~\ref{binary}. Conditions 1 and 2 are clearly 
satisfied and we have to check condition 3.

\begin{tvrz} 
Let $G=([s],E)$ be a simple graph without isolated vertices and $B=\{e_1,e_2,\dots,e_n\}_{\le_l}$ be a subset of $E$. There 
exists a unique simple graph $H=([r],F)$, $F=\{f_1,f_2,\dots,f_n\}_{\le_l}$, of size $n$ without isolated vertices
such that 
$G|\{e_i,e_j\}=F|\{f_i,f_j\}$ for every $1\le i<j\le n$. Moreover, there is an increasing mapping $m:\;[r]\to[s]$ 
such that $m(f_i)=e_i$ for every $1\le i\le n$. 
\end{tvrz}
\duk
$H$ is obtained from $B$ by relabeling the vertices in $V=\bigcup_{e\in B}e$, $|V|=r$, using the unique increasing mapping from $V$
to $[r]$. To construct the mapping $m$, we take the unique $\le_l$-increasing mapping $M:\;F\to E$ 
sending $F$ to $B$ and for a vertex $x\in[r]$ we take an arbitrary edge $f\in F$ with $x\in f$ (since $x$ is not isolated, $f$ 
exists) and define $m(x)=\min M(f)$ if $x=\min f$ and $m(x)=\max M(f)$ if $x=\max f$. Since $M$ preserves types of pairs of 
edges, the value $m(x)$ does not depend on the selection of $f$. Also, $m$ sends $f_i$ to $e_i$ and 
is increasing. The image of each such mapping $m$ is $\bigcup_{e\in B}e$ and $m$ is unique. If $H'$ is another graph with 
the stated property and $m'$ is the corresponding mapping, $m\circ (m')^{-1}$ and $m'\circ m^{-1}$ give an ordered isomorphism 
between $H$ and $H'$. Thus $H$ is unique. 
\kduk

\noindent
We see that simple ordered graphs without isolated vertices, with the ordered subgraph relation and with size being measured by 
the number of edges, form a binary class of objects. $(P,\le_P)$ is the discrete ordering $D_6$.

\bigskip\noindent
{\bf Example 8.  Ordered multigraphs counted by edges.} Let ${\cal O}$ be the set of multigraphs with the vertex set $[n]$ and 
without isolated vertices. The containment $\preceq$ is the ordered subgraph relation and size is the number of edges 
counted with multiplicity. More precisely, in $G=([m], E)\in{\cal O}$ 
we interpret $E$ as a (multiplicity) mapping $E:\;{[m]\choose 2}\to\N_0$, and we have 
$G=([m], E)\preceq H=([n],F)$ iff there is an increasing mapping $f:\;[m]\to[n]$ and an ${m\choose 2}$-tuple 
$\{f_e:\;e\in{[m]\choose 2}\}$ of increasing mappings $f_e:\;[E(e)]\to\N$ such that, for every $e\in{[m]\choose 2}$, 
the image of $f_e$ is a subset of $[F(f(e))]$. The embeddings are the pairs $(f,\{f_e:\;e\in{[m]\choose 2}\})$ and 
$\circ$ is composition of mappings, applied to $f$ and to the mappings $f_e$. 
There is one atom $([2],E)$, where $E([2])=1$, and the size of 
$G=([m], E)$ is the total multiplicity $\sum_{e\subset[m],|e|=2}E(e)$. ${\cal O}_2$ has seven elements. 
The set of atoms $A(G)$ of $G=([m], E)$ can be identified with $\{(e,i):\;e\in{[m]\choose 2},i\in[E(e)]\}$ and the linear 
ordering $(A(G),\le_G)$ is given by $(e,i)\le_G(e',i')$ iff $e<_l e'$ or ($e=e'\;\&\;i\le i'$). The conditions in
Definitions~\ref{class} and \ref{binary} are verified as in the previous example. Therefore multigraphs form a 
binary class of objects. $(P,\le_P)$ is the discrete ordering $D_7$.

\bigskip\noindent
{\bf Example 9. Ordered $k$-uniform hypergraphs counted by edges.} For $k\ge 2$, we generalize
Example 7 to $k$-uniform simple hypergraphs $H=([m],E)$ (so $E\subset{[m]\choose k}$) without isolated vertices. 
The containment $\preceq$ is the ordered subhypergraph relation and size is the number of edges. There is one atom 
$([k],\{[k]\})$. It is not hard to count that ${\cal O}_2$ has 
$$
r=r(k)=\sum_{m=0}^{k-1}{k-1\choose m}\left({2k-m-1\choose k-1}+\frac{1}{2}{2k-m-2\choose k-1}\right)-\frac{1}{2}
$$
elements. $(P,\le_P)$ is the discrete ordering $D_r$.

\bigskip\noindent
{\bf Example 10. Words with the subsequence relation.} For a finite alphabet $A$, let ${\cal O}$ be the set of all
words $u=a_1a_2\dots a_n$ over $A$ and $\preceq$ be the subsequence relation, 
$a_1a_2\dots a_m\preceq b_1b_2\dots b_n$ iff there exists an $m$-tuple $1\le j_1<j_2<\dots<j_m\le n$ 
such that $a_i=b_{j_i}$ for $1\le i\le m$. There are $|A|$ atoms and ${\cal O}_2$ has $r=|A|^2$ elements.   
$(P,\le_P)$ is the discrete ordering $D_r$.

\section{The constant and Fibonacci dichotomies for binary classes of objects}

In this section we prove for $({\cal C}(P),\preceq)$ in Theorem \ref{bounded} the constant dichotomy and in 
Theorem \ref{fibonacci} the Fibonacci dichotomy. Both proofs can be read independently. 
$P$ denotes a finite $l$-element poset on $[l]$ and $l$ is 
always the number of colors. We work with the class $({\cal C}(P),\preceq)$ of all edge 
$P$-colored complete graphs $(n,\chi)$, $n\in\N_0$ and $\chi:\;{[n]\choose 2}\to [l]$. Recall that
$${\textstyle
(m,\psi)\preceq(n,\chi)\iff \exists\mbox{ increasing }f:\;[m]\to[n],\ 
\psi(e)\le_P \chi(f(e))\;\forall e\in{[m]\choose 2}.}
$$
Let $K=(n,\chi)$ be a coloring. The {\em reversal} of $K$ is the coloring $(n,\psi)$ where 
$\psi(\{i,j\})=\chi(\{n-i+1,n-j+1\})$. If $A\subset[n]$ and $\chi|{A\choose 2}$ is constant, we call $A$ a ($\chi$-) 
{\em homogeneous} or ($\chi$-) {\em monochromatic} set. We denote by $R(a;l)$ the Ramsey number for pairs and $l$ colors;  
$R(a;l)$ is the smallest $n\in\N$ such that every coloring $\chi:\;{[n]\choose 2}\to[l]$ has a $\chi$-homogeneous 
set $A\subset[n]$ with size $|A|=a$ (Ramsey \cite{rams}, Graham, Spencer and Rothschild \cite{grah_spen_roth}, 
Ne\v set\v ril \cite{nese}). 

\begin{veta}\label{bounded}
If $X$ is an ideal in $({\cal C}(P),\preceq)$ then either $|X_n|$ is constant for all $n>n_0$ or $|X_n|\ge n$ 
for all $n\ge 1$. 
\end{veta}

\noindent
By Proposition~\ref{discord}, it suffices to prove this if $P$ is a discrete ordering $D_P$. 
We cannot use Corollary~\ref{dvebarvy} to reduce the situation to two colors because we want to prove a result stronger than $|X_n|\ll 1$ but the argument for $l$ colors 
is not too much harder than for two. We need some definitions and auxiliary results.

We say that a coloring $(n,\chi)$ is $r$-{\em rich}, where $r\ge 1$ is an integer, 
if $n=2r-1$ and one the following 
two conditions holds. In {\em type 1} $r$-rich coloring, in $(n,\chi)$ or in its reversal we have $\chi(\{i,i+1\})=a$ for 
$1\le i\le r-1$ and $\chi(\{r,r+1\})=b$ for two colors $a\ne b$. In {\em type 2} $r$-rich coloring, in $(n,\chi)$ or in 
its reversal we have  $\chi(\{1,i\})=a$ for $2\le i\le r$ and $\chi(\{1,r+1\})=b$ for two colors $a\ne b$. We impose no 
restriction on colors of the remaining ${n\choose 2}-r$ edges.

\begin{lema}\label{rich}
If the ideal $X$ contains for every $r\ge 1$ an $r$-rich coloring then $|X_n|\ge n$ for all $n\ge 1$.
\end{lema}
\duk 
If $K=(n,\chi)\in X$ is $r$-rich of type 1, the $r$ restrictions of $K$ to $[i,i+r-1]$ (or to $[n-i-r+2,n-i+1]$ if $K$ is 
reversed) for $1\le i\le r$ are mutually distinct and show that $|X_r|\ge r$. The argument for type 2 $r$-rich colorings 
is similar.
\kduk

\noindent
Note that the containment of an $r$-rich coloring for all $r\ge 1$ is equivalent with the containment for infinitely many $r$
because every $r$-rich coloring contains an $s$-rich coloring for $s=1,2,\dots,r$.

We say that a coloring $(n,\chi)$ is $r$-{\em simple}, where $r\ge 1$ is an integer, if $[r+1,n-r]$ is 
$\chi$-homogeneous and for every fixed $v\in[r]\cup[n-r+1,n]$ the $n-4r$ edges $\{v,w\}$, 
$w\in[2r+1,n-2r]$, have in 
$\chi$ the same color. By the definition, every coloring $(n,\chi)$ with $n\le 2r+2$ is $r$-simple.
We say that a set $X$ of colorings is $r$-simple if each coloring in $X$ is $r$-simple. 

\begin{lema}\label{simple}
If an ideal of colorings $X$ is $r$-simple then there is a constant $d\in\N_0$ such that $|X_n|=d$ for all $n>n_0$.
\end{lema}
\duk
Colorings which are $r$-simple enjoy this property: If $n\ge 4r+2$ and  
$(n,\chi_1)$ and $(n,\chi_2)$ are two distinct $r$-simple colorings, then their restrictions to 
$[n]\backslash\{2r+1\}$ are also distinct. Thus for all $n\ge 4r+2$ the restrictions of the colorings in $X_n$ to 
$[n]\backslash\{2r+1\}$ are mutually distinct and show that $|X_{n-1}|\ge |X_n|$, which implies the claim.
\kduk

\noindent
Theorem~\ref{bounded} now follows from the next proposition.

\begin{tvrz}\label{tezsi}
For every $r\in\N$ there is a constant $c=c(r)\in\N$ such that every ideal of colorings contains an $r$-rich 
coloring or is $c$-simple.
\end{tvrz}

\noindent
{\bf Proof of Theorem~\ref{bounded}.} Let $X$ be an ideal in ${\cal C}(P)$. If $X$ contains an $r$-rich coloring for 
every $r\ge 1$, then $|X_n|\ge n$ for every $n$ by Lemma~\ref{rich}. If not, then by Proposition~\ref{tezsi} $X$ is $c$-simple 
for some $c$ and by Lemma~\ref{simple} $|X_n|$ is constant from some $n$ on. 
\kduk

\noindent
For the proof of Proposition~\ref{tezsi} we shall need three lemmas on situations forcing appearance of $r$-rich colorings.

\begin{lema}\label{homoginterval}
Let $r\ge 1$ be an integer, $(n,\chi)$ be a coloring, $A\subset[n]$ be a $\chi$-homogeneous set with 
the maximum cardinality, and $A'\subset A$ be obtained from $A$ by deleting the first $2r-2$ and the last $2r-2$ elements. 
Suppose further that $A'$ is not an interval in $[n]$. Then $(n,\chi)$ contains an $r$-rich coloring. 
\end{lema}
\duk
We denote the set of the first (last) $r-1$ elements of $A$ by $B_1$ ($B_2$) and the set of the first (last) $2r-2$ elements 
of $A$ by $C_1$ ($C_2$). The assumption on $A'=A\backslash(C_1\cup C_2)$ implies that 
there is an $x\in[n]\backslash A$ such that $C_1<x<C_2$. Since $|A|$ is maximum, there is a $y\in A$ such that the color of
$\{x,y\}$ is different from the color of the edges lying in $A$. 
If $y\in B_1$, then $y$, $C_1\backslash B_1$, $x$, and the next $r-2$ elements of 
$A$ after $x$ form (with restricted $\chi$) an $r$-rich coloring of type 2. If $y\not\in B_1$ but $y<x$, then $B_1$, $y$, 
$x$, and the next $r-2$ elements of $A$ after $x$ form an $r$-rich coloring of type 1. The case when $y>x$ is symmetric 
and is treated similarly.
\kduk

\begin{lema}\label{konce}
Let $r\ge 1$ be an integer, $(n,\chi)$ be a coloring, $s$ be the maximum size of a $\chi$-homogeneous subset of $[n]$, 
$A\subset[n]$ be a $\chi$-homogeneous set with $|A|=s-(4r-4)$, $B\subset[n]$ be a $\chi$-homogeneous set with $A<B$, and let 
$|A|\ge 2r,|B|\ge 6r$. Then $(n,\chi)$ contains an $r$-rich coloring. 
\end{lema}
\duk
Let $a$, respectively $b$, be the color of the edges lying in $A$, respectively 
in $B$. If $a\ne b$, then 
the last $r$ elements of $A$ and the first $r-1$ elements of $B$, or the first $r$ elements of $B$ and the last $r-1$ elements of
$A$ form an $r$-rich coloring of type 1 (depending on whether $\chi(\{\max A,\min B\})$ differs from $a$ or from $b$).

Suppose that $a=b$. Since $|A|+|B|>s$, $A\cup B$ is not homogeneous and there exist 
$x\in A$ and $y\in B$ such that 
$\chi(\{x,y\})\ne a$. Let $y-x$ be minimum, that is, if $x\le x'\in A$, $B\ni y'\le y$, and at least 
one inequality is strict, then $\chi(\{x',y'\})=a$. We show that any position of $x$ and $y$ produces an $r$-rich coloring.
We denote by $C_1$ ($C_2$) the set of the first (last) $r-2$ elements of $A$ ($B$).
Suppose first that $x\not\in C_1$. If $y$ is among the last $r-1$ elements of $B$, then
$y$, the previous $r-1$ elements of $B$, $x$, and $C_1$ form an $r$-rich coloring of type 2.
If $y$ is not among the last $r-1$ elements of $B$, then these elements, $y$, $x$, and $C_1$ 
form an $r$-rich coloring of type 1. A symmetric argument shows that if $y\not\in C_2$ then 
we have an $r$-rich coloring. The remaining case when $x\in C_1$ and $y\in C_2$ does not occur because then, 
by the minimality of the length of $\{x,y\}$, $(A\backslash C_1)\cup(B\backslash C_2)$ would be a homogeneous set with 
size $|A|+|B|-2(r-2)\ge |A|+4r+4=s+8$, contradicting the definition of $s$.
\kduk

\begin{lema}\label{treti}
Let $r\ge 1$ be an integer, $(n,\chi)$ be a coloring, $v\in [n]$, $A\subset[n]$ be a set such that $v<A$ or $v>A$, and $B\subset A$
be obtained from $A$ by the deletion of the first $l(r-2)+1$ and the last $l(r-2)+1$  elements. Suppose further that not all edges 
$\{v,w\}$, $w\in B$, have in $\chi$ the same color. Then $(n,\chi)$ contains an $r$-rich coloring.
\end{lema}
\duk
Let $v<A$, the proof for $v>A$ is very similar. By the assumption there is a $w\in B$ such that 
$b=\chi(\{v,w\})\ne a=\chi(\{v,\max B\})$. By the pigeonhole principle, 
some $r-1$ edges $\{v,z_1\},\dots,\{v,z_{r-1}\}$, where $z_i\in A$ and $z_i<B$, have the same color $c$. Thus 
$v$, $z_1,\dots,z_{r-1}$, $w$ or $\max B$ (depending on whether $c\ne b$ or $c\ne a$), and the last $r-2$ elements of
$A$ form an $r$-rich coloring of type 2. 
\kduk

\bigskip\noindent
{\bf Proof of Proposition~\ref{tezsi}.} 
We assume that $X$ is an ideal in ${\cal C}(P)$ which contains no $r$-rich coloring for some $r\ge 2$. We show that then
$X$ must be $c$-simple for
$$
c=\max(R(6r;l),l(r-2)+1)
$$
where $R(\cdot;\cdot)$ is the Ramsey number.
Let $(n,\chi)\in X$ be arbitrary. We may assume that $n>2c+2$. We take a set $A\subset[n]$ obtained from a $\chi$-homogeneous 
subset with the maximum cardinality by deleting the first $2r-2$ and the last $2r-2$ elements. By Lemma~\ref{homoginterval},  
the Ramsey theorem, and Lemma~\ref{konce}, $A$ is an interval in $[n]$ and $\min A<c+1$, $\max A>n-c$. 
Thus $[c+1,n-c]$ is $\chi$-homogeneous. Let 
$v\in[c]\cup[n-c+1,n]$ be arbitrary. By Lemma~\ref{treti} applied to $v$ and $[c+1,n-c]$, all edges $\{v,w\}$, $w\in[2c+1,n-2c]$, 
must have the same color. Thus $(n,\chi)$ is $c$-simple.
\kduk

\noindent
This completes the proof of Theorem~\ref{bounded}.

\begin{veta}\label{fibonacci}
If $X$ is a ideal in $({\cal C}(P),\preceq)$ then either $|X_n|\le n^c$ for all $n\ge 1$ for a constant $c>0$, or
$|X_n|\ge F_n$ for all $n\ge 1$ where $(F_n)_{n\ge 1}=(1,2,3,5,8,13,\dots)$ are the Fibonacci numbers. 
\end{veta}

\noindent
As in the proof of Theorem~\ref{bounded}, we define ``wealthy''
colorings (of four types) and ``tame'' colorings and show that colorings with unbounded wealth produce growth 
at least $F_n$ and that bounded tameness admits only polynomially many colorings. The proof is completed by 
showing that the colorings in any ideal are either unboundedly wealthy or boundedly tame.
By Corollary~\ref{dvebarvy} and the following remark, it suffices to prove the theorem only for the two-element discrete poset 
$P=D_2$, that is, for graphs and ordered induced subgraph relation. To make explicit the symmetry between edges and 
nonedges in this case, we prefer to use the language of colorings. Therefore by a {\em coloring} we shall mean in the proof always a black-white 
edge coloring $(n,\chi)$ of a complete graph, $\chi:\;{[n]\choose 2}\to\{black,white\}$, and if we use two 
distinct colors $c,d$, one should bear in mind that $\{c,d\}=\{black,white\}$. 

Let $r\in\N$. A coloring $K=(r,\chi)$ is $r$-{\em wealthy of type 1} if in $K$ or in its reversal we have 
$\chi(\{1,i\})\ne\chi(\{1,i+1\})$ for all $i\in[2,r-1]$. 
$K=(3r,\chi)$ is $r$-{\em wealthy of type 2 } if none of the $r$ consecutive triangles 
$\{3i-2,3i-1,3i\}$, $1\le i\le r$, is $\chi$-homogeneous. We use two incarnations of the Fibonacci number $F_n$. 

\begin{lema}\label{fibo}
(i) $F_n$ equals to the number of 0-1 strings $s_1s_2\dots s_{n-1}$ with no two consecutive 1s, i.e., 
avoiding the pattern $s_is_{i+1}=11$. (ii) $F_n$ equals to the number of 0-1 strings $s_1s_2\dots s_{n-1}$ 
avoiding the patterns $s_{2i-1}s_{2i}=01$ and $s_{2i}s_{2i+1}=10$.
\end{lema}
\duk
Both results are easily proved by induction on $n$. 
\kduk

\noindent
We call the strings in (i) {\em fib1 strings} and the strings in (ii) {\em fib2 strings}.

\begin{lema}\label{fibobound12}
If there is an $i\in\{1,2\}$ so that the ideal $X$ contains for every $r\ge 1$ an $r$-wealthy coloring of type $i$, 
then $|X_n|\ge F_n$ for all $n\ge 1$.
\end{lema}
\duk
If $X$ contains for every $r\ge 1$ an $r$-wealthy coloring of type 1, it follows that for every $n\in\N$ and every
subset $A\subset [2,n]$ there exists a coloring $K_A=(n,\chi)$ in $X$ such that $\chi(\{1,i\})=black\iff i\in A$, 
or the same holds for the reversals of $K_A$s. Because for fixed $n$ all $2^{n-1}$ colorings $K_A$ are mutually distinct, we 
have $|X_n|\ge 2^{n-1}\ge F_n$. 

Suppose that $X$ contains for every $r\ge 1$ an $r$-wealthy coloring of type 2. Using the pigeonhole principle and 
the Ramsey theorem, we regularize the situation and obtain the colorings $(3,\phi)$, $(6,\psi)$, and $(3r,\chi_r)$, 
$r=1,2,\dots$, which all lie in $X$ and are such that in $(3r,\chi_r)$ all triangles $T_i=\{3i-2,3i-1,3i\}$, 
$1\le i\le r$, and the edges between them are colored in the same way and independently of $r$: if $\{a,b\}\subset T_i$ then 
$\chi_r(\{a,b\})=\phi(\{a-3(i-1),b-3(i-1)\})$ and if $a\in T_i$ and $b\in T_j$, $1\le i<j\le r$, then 
$\chi_r(\{a,b\})=\psi(\{a-3(i-1),b-3(j-2)\})$. The coloring $(3,\phi)$ of the triangles is not monochromatic and thus 
there is an edge $\{a,b\}\subset T_1$ such that not all of the four edges connecting $\{a,b\}$ and $\{a+3,b+3\}$ have 
color $c=\phi(\{a,b\})$. It follows that there are colorings $(4,\kappa)$ and $(2r,\lambda_r)$, $r=1,2,\dots,$ 
which lie in $X$ and are such that 
(i) the edges $\{1,2\},\{3,4\},\dots,\{2r-1,2r\}$ have in $\lambda_r$ the same color, say black, (ii) if $a\in\{2i-1,2i\}$ and 
$b\in\{2j-1,2j\}$ with $1\le i<j\le r$, then $\lambda_r(\{a,b\})=\kappa(\{a-2(i-1),b-2(j-2)\})$, and (iii) at least one 
of the four edges $\{1,3\},\{1,4\},\{2,3\},\{2,4\}$ is in $\kappa$ colored
white. Suppose, for example, that $\{1,4\}$ is white. If $\{1,3\}$ is black, it follows that $(2r,\lambda_r)$ 
contains an $r$-wealthy coloring of type 1 and we are in the previous case. This argument shows that we may assume that 
in $(2r,\lambda_r)$ all edges $\{1,2\},\{3,4\},\dots,\{2r-1,2r\}$ are black and all other edges white. 
It follows that for every $n\in\N$ and every fib1 string $w=s_1s_2\dots,s_{n-1}$ there is a coloring 
$K_w=(n,\chi)\in X$ such that $\chi(\{i,i+1\})=black\iff s_i=1$. Since for distinct $w$s the corresponding
colorings $K_w$ are distinct as well, by (i) of Lemma~\ref{fibo} we conclude that $|X_n|\ge F_n$. 
\kduk

Before defining wealthy colorings of types 3 and 4, we introduce notation on 0-1 matrices which we will use
to represent colorings.
If $M$ is an $r\times s$ 0-1 matrix, any row and column of $M$ consists of alternating intervals of consecutive 
0s and 1s. Let $al(M)$ be the maximum number of these intervals in a row or in a column, taken over all $r+s$ rows and columns. For every
$j\in[s]$ we let $C(M,j)\subset[r]$ be the row indices of the lowest entries of these intervals in the $j$-th 
column, with $r$ omitted: $a\in C(M,j)$ iff $M(a,j)\ne M(a+1,j)$. We denote $C(M)=\bigcup_{j=1}^s C(M,j)$. For a 
coloring $K=(2r,\chi)$ we define the $r\times r$ 0-1 matrix $M_K$ by $M_K(i,j)=0$ iff $\chi(\{i,r+j\})=white$. 
Similarly, if $K=(n,\chi)$ is a coloring and $I=\{x_1<x_2<\dots<x_r\}<J=\{y_1<y_2<\dots<y_s\}$ 
are two subsets of $[n]$, we define the $r\times s$ 0-1 matrix $M_{I,J}$ by $M_{I,J}(i,j)=0$ iff 
$\chi(\{x_i,y_j\})=white$; we suppress in notation the dependence on $K$ which will be clear from the context.

We say that $C=(c_1,c_2,\dots,c_k)$, with $c_i\in[r]\times[s]$ being in the (row,column) coordinates format, is a 
{\em southeast path} in an $r\times s$ 0-1 matrix $M$
if in $C$ alternate south and east steps and $C$ starts with a south step: $c_{2i}-c_{2i-1}\in\N\times\{0\}$ and 
$c_{2i+1}-c_{2i}\in\{0\}\times\N$. If $M=M_K$ or $M=M_{I,J}$ for some coloring $K$ then $C$ corresponds to a path 
in the coloring, with the $k$ edges
$$
\{a_1,b_1\}, \{a_2,b_1\}, \{a_2,b_2\}, \{a_3,b_2\}, \{a_3,b_3\},\dots, \{a_p,b_q\}
$$ 
where $p=\lceil (k+1)/2\rceil$, 
$q=\lceil k/2\rceil$, and  $1\le a_1<a_2<\dots<a_p<b_1<b_2<\dots<b_q$. We call such paths 
{\em back-and-forth paths}. 

If $M$ is an $r\times s$ 0-1 matrix and $M'$ is an $r'\times s'$ 0-1 matrix, we say that $M'$ is contained in $M$, 
$M'\preceq M$, if there are increasing injections $f:\;[r']\to[r]$ and $g:\;[s']\to[s]$ such that 
$M(f(i),g(j))=M'(i,j)$ for all $i\in[r'],j\in[s']$. We say that $M'$ is a {\em submatrix} of $M$. 
If $M'\preceq M$ and $M=M_{I,J}$ for two subsets $I<J$ in $[n]$ and 
a coloring $(n,\chi)$, then there are subsets $I'\subset I$ and $J'\subset J$ such that $M'=M_{I',J'}$.
We denote by $I_r$ the $r\times r$ identity matrix with 1s on 
the main diagonal and 0s elsewhere and by $U_r$ the upper triangular $r\times r$ matrix with 1s 
above and on the main diagonal and 0s below it. We call two $r\times s$ 0-1 matrices $M$ and $M'$ {\em similar} if 
$M'=M$ or $M'$ is obtained from $M$ by the vertical mirror image and/or swapping 0 and 1

We say that a coloring $K=(2r,\chi)$ is $r$-{\em wealthy of type 3} if $M_K$ is similar to $I_r$. 
$K=(2r,\chi)$ is $r$-{\em wealthy of type 4} if $M_K$ is similar to $U_r$. Note that for $i\in\{1,2,3,4\}$ and 
$r\in\N$, every $r$-wealthy coloring of type $i$ contains an $s$-wealthy coloring of type 
$i$ for $s=1,2,\dots,r$ and so for an ideal $X$ to contain an $r$-wealthy coloring of type $i$ for every $r\ge 1$ is equivalent with containing it for infinitely many $r$. 

\begin{lema}\label{fibobound34}
If there is an $i\in\{3,4\}$ so that the  set $X$ contains for every $r\ge 1$ an $r$-wealthy coloring of type $i$, 
then $|X_n|\ge F_n$ for all $n\ge 1$.
\end{lema}
\duk
Let $X$ contain for every $r\ge 1$ an $r$-wealthy coloring $K_r$ of type 3. We may assume that always $M_{K_r}=I_r$. 
It can be seen that if $n\in\N$ and $w=s_1s_2\dots s_{n-1}$ is any fib1 string, then 
for $r\ge 2n$ one can draw in $I_r$ a southeast path $(c_1,c_2,\dots,c_{n-1})$ such that 
$I_r(c_i)=s_i$. Thus for every $w$ there is a coloring $K_w=(n,\chi)\in X$ whose unique spanning 
back-and-forth path is colored according to $w$. By (i) of Lemma~\ref{fibo} we have $|X_n|\ge F_n$. 

Let $X$ contain for every $r\ge 1$ an $r$-wealthy coloring $K_r$ of type 4. We may assume that always 
$M_{K_r}=U_r$. It can be seen that if $n\in\N$ and $w=s_1s_2\dots s_{n-1}$ is any fib2 string, then 
for $r\ge 2n$ one can draw in $U_r$ a southeast path $(c_1,c_2,\dots,c_{n-1})$ such that 
$U_r(c_i)=s_i$. Again, by (ii) of Lemma~\ref{fibo} we have $|X_n|\ge F_n$. 
\kduk

\begin{lema}\label{alternation}
Let $M$ be an $r\times s$ 0-1 matrix that satisfies $al(M)\le k$ and $|C(M)|\le l$, and $a$ be the 
number of the column indices $j\in[s]$ for which the $j$-th column of $M$ differs from the $(j+1)$-th one. 
Then
$$
a\le (k-1)(2l+1). 
$$
\end{lema}
\duk
The $j$-th column of $M$ is uniquely determined by $C(M,j)\subset C(M)$ and 
by $M(1,j)\in\{0,1\}$. It follows that any two different columns in $M$ 
must differ in entries with row index lying in the set 
$D=C(M)\cup\{i+1:\;i\in C(M)\}\cup\{1\}$, which has at most $2l+1$ elements. 
By the pigeonhole principle, if $a>(k-1)(2l+1)$ then there are $k$ column indices 
$1\le j_1<j_2<\dots<j_k<s$ and a row index $b\in D$ such that 
$M(b,j_i)\ne M(b,j_i+1)$ for all $i\in[k]$. Thus the $b$-th row of $M$ 
consists of at least $k+1$ alternating intervals 
of 0s and 1s, which contradicts $al(M)\le k$.
\kduk

\begin{lema}\label{matrices}
Let $(M_n)_{n\ge 1}$ be an infinite sequence of 0-1 matrices such that (i) the sequence $(al(M_n))_{n\ge 1}$ is bounded but (ii) 
$(|C(M_n)|)_{n\ge 1}$ is unbounded. Then either (a) for every $r$ there is an $n$ and a matrix $I_r'$ similar to $I_r$ 
such that $I_r'\preceq M_n$ or (b) for every $r$ there is an $n$ and a matrix $U_r'$ similar to $U_r$ such that $U_r'\preceq M_n$.
\end{lema}
\duk
We prove the result under the weaker assumption with $al(M_n)$ replaced by $al_c(M_n)$ that is defined by taking 
the maximum (of the numbers of intervals of consecutive 0s and 1s) only over the columns of $M_n$. Using the pigeonhole 
principle and replacing $(M_n)_{n\ge 1}$ by an appropriate subsequence of submatrices, we may assume in addition to 
(ii) that there is an $s\ge 1$ and a $c\in\{0,1\}$ such that the first row of every $M_n$ contains only $c$s and
$|C(M_n,j)|=s$ for every $n\ge 1$ and $j$. We set $C(M_n,j)=\{r_{n,j,1}<r_{n,j,2}<\dots<r_{n,j,s}\}$ and denote $c_n$ 
the number of columns in $M_n$. We proceed by induction on $s$. It is clear that if $s=1$ then the sequence
$S=(|\{r_{n,j,1}:\;1\le j\le c_n\}|)_{n\ge 1}$ is unbounded. Suppose that $s\ge 2$ and $S$ is bounded. Taking
a subsequence of submatrices, we may then assume in addition to (ii) that $r_{n,j,1}=r_n$ for 
$1\le j\le c_n$ and all $n$. We take from every $M_n$ only rows $r_n+1, r_n+2,\dots$ and obtain a sequence of
matrices $(N_n)_{n\ge 1}$ satisfying $|C(N_n,j)|=s-1$ for every $n,j$ and (ii), which means that we are done by induction. Thus we may
assume that $S$ is unbounded even for $s\ge 2$. We take a subsequence of submatrices once more and may assume that $(M_n)_{n\ge 1}$ 
satisfies: $|C(M_n,j)|=s$ for all $n$ and $j$, $(c_n)_{n\ge 1}$ is unbounded, and for every $n$ the $c_n$ row indices $r_{n,j,1}$, 
$j\in[c_n]$, are mutually distinct.

Suppose that $s=1$. Using the Erd\H os-Szekeres lemma, we may assume that moreover
for every $n$ the sequence $(r_{n,j,1}:\;j=1,2,\dots,c_n)$ is strictly increasing or that for every $n$ it is strictly decreasing. 
Keeping from $M_n$ only the rows $r_{n,1,1}, r_{n,2,1},\dots,r_{n,c_n,1}$ (and all columns), 
we obtain a matrix similar to $U_{c_n}$. We see 
that (b) holds. In the case that $s\ge 2$ we denote by $I_{n,j}$ the interval $[r_{n,j,1}+1,r_{n,j,s}]$ and by $i(n)$
(resp. $d(n)$) the maximum number of intervals among $I_{n,1},I_{n,2},\dots,I_{n,c_n}$ which share one point (resp. which are 
mutually disjoint). It follows that $(i(n))_{n\ge 1}$ or $(d(n))_{n\ge 1}$ is unbounded. In the former case we may assume, 
turning to a subsequence of submatrices, that $r_n\in\bigcap_{j=1}^{c_n}I_{n,j}$ for every $n\ge 1$ for some row indices $r_n$.
Keeping from $M_n$ only the rows $1,2,\dots,r_n$, we obtain a sequence of matrices $(N_n)_{n\ge 1}$ which satisfies 
$al_c(N_n)\le s-1$ for every $n$ and (ii), which means that we are done by induction. In the latter case we may assume, using
the Erd\H os-Szekeres lemma and turning to a subsequence of submatrices, that for every $n$ we have 
$I_{n,1}<I_{n,2}<\dots<I_{n,c_n}$ or that for every $n$ we have $I_{n,1}>I_{n,2}>\dots>I_{n,c_n}$. We select row indices 
$t_{n,j}\in I_{n,j}$ such that, for all $n$ and $j\in[c_n]$, $M_n(t_{n,j},j)\ne M_n(r_{n,j,1},j)=M_n(r_{n,j,s}+1,j)$ if $s$ is even
and $M_n(t_{n,j},j)=M_n(r_{n,j,1},j)\ne M_n(r_{n,j,s}+1,j)$ if $s$ is odd. Keeping from $M_n$ only the $c_n$ rows 
$t_{n,1},t_{n,2},\dots,t_{n,c_n}$, we obtain a matrix similar to $I_{c_n}$ (for even $s$) or to $U_{c_n}$ (for odd $s$). Thus 
(a) or (b) holds.
\kduk

A coloring $(n,\chi)$ is $m$-{\em tame}, where $m\in\N$, if the following three conditions are 
satisfied.
\begin{enumerate}
\item There is an interval partition $I_1<I_2<\dots<I_s$ of $[n]$ such that $s\le m$ and every $I_i$ is 
$\chi$-monochromatic.
\item For every two subintervals $I<J$ in $[n]$ we have $al(M_{I,J})\le m$.
\item For every two subintervals $I<J$ in $[n]$ we have $|C(M_{I,J})|\le m$.
\end{enumerate}
A set of colorings $X$ is $m$-tame if every coloring $(n,\chi)\in X$ is $m$-tame.

\begin{lema}\label{polybound}
For every $m\in\N$ there is a constant $c=c(m)$ such that the number of $m$-tame colorings $(n,\chi)$ is bounded by $n^c$.
\end{lema}
\duk
The partition $I_1<I_2<\dots<I_s$ of $[n]$ satisfying condition 1 and the $s$ colors $\chi|{I_i\choose 2}$ 
can be selected in $c_1=\sum_{s=1}^m{n-1\choose s-1}2^s\le (2n)^m$ ways. The colors of the remaining edges in $(n,\chi)$
are determined by the 0-1 matrices $M=M_{I_u,I_v}$, $1\le u<v\le s\le m$. Let us bound, for fixed $u,v$, the number of 
matrices $M$ satisfying conditions 2 and 3. The number of possibilities for one column of $M$ is by condition 2 at most 
$c_2=2\sum_{i=1}^m{p-1\choose i-1}\le (2n)^m$, where $p=|I_u|\le n$, and all $q$ columns of $M$, $q=|I_v|\le n$, can be 
selected by Lemma~\ref{alternation} in at most 
$$
c_3=\sum_{i=1}^{2m^2}
{q-1\choose i-1}c_2^i\le 2m^2\cdot(n-1)^{2m^2-1}\cdot(2n)^{2m^3}\le
(2n)^{4m^3}
$$ 
ways. The total number of $m$-tame colorings $(n,\chi)$ is therefore at most
$$
c_1c_3^{{s\choose 2}}\le c_1c_3^{{m\choose 2}}\le (2n)^{2m^5+m}\le n^c,\ n\ge 2.
$$
\kduk

Let $K=(n,\chi)$ be a coloring. The {\em interval decomposition} of $K$ is the unique partition of $[n]$ 
in nonempty intervals $I_1<I_2<\dots<I_s$ defined as follows. $I_1$ is the longest initial interval such that $I_1$ is 
$\chi$-monochromatic, $I_2$ is the longest following interval such that $I_2$ is $\chi$-monochromatic, and so on. Clearly $|I_i|\ge 2$ 
for all $i<s$. We let $I(K)=s$ denote the number of intervals in the decomposition.

\begin{tvrz}\label{notmtame}
If $X$ is an ideal in $({\cal C}(D_2),\preceq)$ that is not $m$-tame for any 
$m$, then there is an 
$i\in\{1,2,3,4\}$ such that $X$ contains for every $r\ge 1$ an $r$-wealthy coloring of type $i$.
\end{tvrz}
\duk
Suppose there is no $m\in\N$ such that $X$ is $m$-tame. Thus one of the three conditions 
in the definition of tameness is violated for infinitely many $m$ on some colorings in $X$. 
If it is condition 1, the quantity $I(K)$, $K\in X$, is unbounded 
and for every $r\ge 1$ there is a coloring $(n,\chi)\in X$ whose interval decomposition $I_1<I_2<\dots<I_s$ satisfies $s\ge r$. 
By the definition,
for every $i$, $1<i\le s$, there is an $x_{i-1}\in I_{i-1}$ such that $\chi(\{x_{i-1},\min I_i\})$ differs from the color 
$\chi|{I_{i-1}\choose 2}$. Hence the triangles $\{x_{2i-1},y_{2i-1},\min I_{2i}\}$, where $1\le i\le r/2$ and 
$y_{2i-1}\in I_{2i-1}\backslash\{x_{2i-1}\}$ is selected arbitrarily, are not monochromatic in $(n,\chi)$ and $X$ contains for every 
$r\ge 1$ an $r$-wealthy coloring of type 2. If condition 2 is violated infinitely many times, it is easy to see 
that $X$ contains for every $r\ge 1$ an $r$-wealthy coloring of type 1.

We are left with the case when conditions 1 and 2 of tameness are satisfied for all colorings in $X$ with a 
constant $m_0$ but condition 3 is violated for all $m$. This implies that for $n=1,2,\dots$ there are colorings 
$(n,\chi_n)\in X$ and subintervals $I_n<J_n$ in $[n]$ such that the sequence of 0-1 matrices $(M_{I_n,J_n})_{n\ge 1}$
satisfies the hypothesis of Lemma~\ref{matrices}. By the conclusion of the lemma, there is an $i\in\{3,4\}$ such that 
$X$ contains for every $r\ge 1$ an $r$-wealthy coloring of type $i$.
\kduk

\noindent
{\bf Proof of Theorem~\ref{fibonacci}.} If $X$ is an ideal in ${\cal C}(D_2)$ that is $m$-tame for an $m$, then by 
Lemma~\ref{polybound} we have $|X_n|\le n^c$ for all $n\ge 1$ with a constant $c>0$. If $X$ is not $m$-tame for any $m$, 
by Proposition~\ref{notmtame} there is an $i\in\{1,2,3,4\}$ so that $X$ contains for every $r\ge 1$ an 
$r$-wealthy coloring of type $i$. By Lemmas~\ref{fibobound12} and \ref{fibobound34} this means
that for all $n\ge 1$ we have $|X_n|\ge F_n$.
\kduk

\section{Concluding remarks}

We conclude with mentioning a few open problems on growths of ideals of permutations and graphs.

The Stanley-Wilf conjecture (B\'ona \cite{bona97jcta, bona97, bona_book}) asserted that for every permutation 
$\pi$ the number of $n$-permutations not containing $\pi$ is exponentially bounded. Equivalently stated, 
for every ideal of permutations $X$ different from the set of all permutations ${\cal S}$ we have $|X_n|<c^n$ for all $n\ge 1$. The conjecture was proved by Marcus and Tardos in \cite{marc_tard} and therefore now we know that
$$
c(X)=\limsup_{n\to\infty}|X_n|^{1/n}<\infty
$$
for every ideal $X\ne{\cal S}$. However, many interesting and challenging problems on growth of ideals of permutations remain open. The following problem was posed by V. Vatter \cite{elde_vatt}.

\bigskip\noindent
{\bf Problem 1.} Is it true that $\lim_{n\to\infty}|X_n|^{1/n}$ always exists? 

\bigskip\noindent
It was proved by Arratia \cite{arra} in the case $X={\rm Forb}(\{\pi\})$. By the Fibonacci hierarchy 
4 (Introduction), it is also true when $c(X)\le 2$.

\bigskip\noindent
{\bf Problem 2.} What are the constants of growth
$$
C=\{c(X):\;X\mbox{ is an ideal of permutations}\}?
$$
Are all of them algebraic? 

\bigskip\noindent
Similar problem was posed \cite[Conjecture 8.9]{balo-boll-morr06a} for hereditary properties of ordered graphs.

It is easy to find ideals of permutations $X$ such that the function 
$n\mapsto|X_n|$ is, respectively, constant $0$, constant $1$, 
$n\mapsto F_{n,k}$ for any fixed $k\ge 2$, and $n\mapsto 2^{n-1}$. Thus, by the Fibonacci hierarchy 4, 
$$
C\cap [0,2]=\{0,1,\alpha_2,\alpha_3,\alpha_4,\dots, 2\}
$$
where $\alpha_2\approx 1.61803$, $\alpha_3\approx 1.83928$, $\alpha_4\approx 1.92756, \dots$ are the limits
$\alpha_k=\lim F_{n,k}^{1/n}$. By the standard results from asymptotic enumeration, $\alpha_k$ is the largest positive real root 
of $x^k-x^{k-1}-x^{k-2}-\cdots-1$. It follows that $\alpha_k\uparrow 2$. It would be interesting to determine further elements of $C$ lying above the first limit point $2$. 

\bigskip\noindent
{\bf Problem 3.} Show that $\min(C\cap (2,\infty))$ exists. What is this number?

\bigskip
A natural question arises if also the remainig two restrictions on growth of hereditary properties of ordered 
graphs proved by Balogh, Bollob\'as and Morris \cite{balo-boll-morr06a}, the polynomial growth 3 and 
the Fibonacci hierarchy 4, can be extended to edge-colored complete graphs with $l\ge 2$ colors. It is not 
too hard to achieve this for the polynomial growth by elaborating the final ``tame'' part of our proof of 
Theorem~\ref{fibonacci}; we hope to say more on this elsewhere. It is plausible to cojecture that the 
proof of the Fibonacci hierarchy in \cite{balo-boll-morr06a} also can be ``upgraded'' from $l=2$ to $l\ge 2$
but this would require more effort. 

Finally, we present an interesting problem on an exponential-factorial jump in growth due to 
Balogh, Bollob\'as and Morris \cite[Conjecture 2]{balo-boll-morr06b}.  

\bigskip\noindent
{\bf Problem 4.} Let $X$ be a hereditary property of ordered graphs. Prove that either $|X_n|<c^n$ for 
all $n\ge 1$ with a constant $c>1$ or 
$$
|X_n|\ge\sum_{k=0}^{\lfloor n/2\rfloor}{n\choose 2k}k!
$$
for all $n\ge 1$.

\bigskip\noindent
They proved \cite[Theorem 4]{balo-boll-morr06b} this jump for the smaller family of monotone properties of 
ordered graphs (and in fact more generally for hypergraphs).

\end{document}